\documentclass[12pt]{article}%
\usepackage{amssymb}
\usepackage{graphicx}
\usepackage{amsmath}
\usepackage{citesort}
\usepackage{amsfonts}%
\providecommand{\U}[1]{\protect\rule{.1in}{.1in}}
\newtheorem{theorem}{Theorem}
\begin{document}

\title{ Thermoacoustic tomography with detectors on an open curve: an efficient
reconstruction algorithm}
\author{Leonid Kunyansky\\Department of Mathematics \\University of Arizona, Tucson.}
\maketitle

\begin{abstract}
Practical applications of thermoacoustic tomography require numerical
inversion of the spherical mean Radon transform with the centers of
integration spheres occupying an open surface. Solution of this problem is
needed (both in 2-D and 3-D) because frequently the region of interest cannot
be completely surrounded by the detectors, as it happens, for example, in
breast imaging. We present an efficient numerical algorithm for solving this
problem in 2-D (similar methods are applicable in the 3-D case). Our method is
based on the numerical approximation of plane waves by certain single layer
potentials related to the acquisition geometry. After the densities of these
potentials have been precomputed, each subsequent image reconstruction has the
complexity of the regular filtration backprojection algorithm for the
classical Radon transform. The peformance of the method is demonstrated in
several numerical examples: one can see that the algorithm produces very
accurate reconstructions if the data are accurate and sufficiently well
sampled, on the other hand, it is sufficiently stable with respect to noise in
the data.

\end{abstract}

\section*{Introduction}

The thermoacoustic tomography (TAT) is based on measurements of acoustic waves
excited in the patient's body by an external electromagnetic pulse
\cite{Kruger,Xu}. The method is gaining popularity among researchers because
it combines the high resolution of the ultrasound tomography with the high
contrast of the images attainable due to the strong variance in the
electromagnetic properties of the body tissues. One of the most promising
applications of this technique is breast imaging, where absorption of the
electromagnetic energy in tumors is several times higher than that in healthy
tissues. Under certain simplifying assumptions (most importantly, that of a
constant speed of sound in the tissue) the corresponding reconstruction
problem can be reduced to the inversion of the spherical mean Radon transform.
In the present paper we propose an efficient numerical algorithm for solving
this problem in 2-D in the case when the detectors lie on an open curve only
partially surrounding the region of interest.

A significant progress has been achieved recently in mathematics of TAT; an
extended discussion and relevant references can be found in reviews
\cite{AKK,FPR,KuchKu}. In particular, several versions of explicit inversions
formulas \cite{Finch,FPR,Kunyansky,Xu-univ,Burgh} and certain series solutions
\cite{Norton1,Norton2,Kunyansky1} have been obtained for measurement schemes
with the detectors lying either on closed surfaces surrounding the region of
interest (ROI), or on certain open but unbounded surfaces, such as a plane or
infinite cylinder.

However, in most applications the ROI is not the whole human body but rather a
part of it --- as it happens, for example, in mammography. In such situations
only some part of the region boundary can be covered by the detectors. This,
in turn, necessitates the development of algorithms that can reconstruct the
image from such data.

The simplest approach to image reconstruction from open surfaces is to place
detectors on a truncated plane or cylinder and to apply the known inversion
formula(s) for the full plane or infinite cylinder. Of course, the
reconstruction would not be exact; moreover, it has been shown that for a
stable reconstruction of a compactly supported function from its spherical
means it is necessary to know the integrals over sufficiently many spheres, so
that the so-called "visibility" condition is satisfied
\cite{LQ,XWAK1,XWAK2,Palambook,Q90}. (It have been shown, in particular, that
under this condition the reconstruction problem can be reduced either to the
inversion of an elliptic operator or, after a certain filtration, to a
solution of the Fredholm integral equation of the second kind
\cite{Q90,Palambook,Palam1,KLM}. In both cases the problem is only mildly
unstable, similarly to the inversion of the standard Radon transform
\cite{Natterer}). The "visibility" condition (for TAT) is satisfied if for
each point $x$ in the ROI each straight line passing through $x$ intersects
the measuring surface at least once. A sphere surrounding the object satisfies
this condition while an infinite plane and an infinite cylinder do not. For a
truncated plane or truncated cylinder the set of "bad" directions is even
larger than for their infinite counterparts. If this condition is violated, it
is practically impossible to accurately reconstruct sharp interfaces (material
boundaries) at those points of the image for which some of the straight lines
normal to the interface do not intersect the measuring surface. (In 2-D an
exact reconstruction technique was developed in \cite{Palam} for the case when
the centers of the integration circles lie on a segment of a straight line.
However, the reconstruction is still, in general, unstable in this case.)

For a given bounded ROI there exist many bounded open acquisition surfaces
that do satisfy the visibility condition. Almost all known reconstruction
techniques applicable to such surfaces are of approximate nature. For example,
by "approximating" the integration spheres by planes and by applying some
version of the classical inverse Radon transform, one can reconstruct an
"approximation" to the image. Due to the symmetry in the classical Radon
projections, the normals to the integration planes should fill only a half of
a unit sphere, making possible the reconstruction from an open measurement
surface. A more sophisticated approach is represented by the so-called
"straightening" methods \cite{PopSush,PopSush2} based on the approximate
reconstruction of the classical Radon projections from the measurements that
correspond to the spherical mean Radon transform. However, these methods do
not yield the exact solution but rather a parametrix of the problem; the image
is reconstructed only up to a certain smooth term. In other words, while jumps
corresponding to sharp material interfaces are reconstructed accurately, the
accuracy of the lower spatial frequencies can not be guaranteed. Unlike the
approximations resulting from the discretization of the exact inversion
formulas (in the situations when such formulas are known), the parametrix
approximations do not converge when the discretization of the data is refined
(in the absence of noise). In \cite{Paltaufalg}, an exact inversion formula
for the spherical surface is used to obtain approximate reconstructions from
other measurement surfaces; in order to further improve the approximation,
additional corrections are introduced. These methods yield different
parametrices than the one proposed in \cite{PopSush,PopSush2}; which one is
better remains an open question.

An accurate numerical reconstruction from an open surface acquisition scheme
was demonstrated in \cite{Anas1} (see also \cite{Anas2}). In this work a
version of an iterative algebraic reconstruction algorithm was successfully
employed to recover a numerically generated phantom. Iterative algebraic
reconstruction algorithms are, however, notoriously slow; the above-mentioned
reconstruction, for example, required the use of a cluster of computers and
took 100 iterations to converge. A faster converging algorithm can be obtained
by combining iterative refining with a parametrix-type algorithm
\cite{Paltauf,Paltaufalg}. These methods are closely related to the general
scheme proposed in \cite{Beylkin} for the inversion of the generalized Radon
transform with integration over general manifolds. It reduces the problem to
the Fredholm integral equation of the second kind, which is well suited for
numerical solution. Such an approach can be viewed as using a parametrix
method as an efficient preconditioner for an iterative solver, which, as it
has been shown in numerical experiments, significantly accelerates the
convergence of the iterations..

Finally, in \cite{Sarah}\ an interesting attempt has been made to generate the
absent data from the consistency conditions on the spherical mean Radon
transform, in order to "numerically close" the open measurement surface. The
resulting algorithm, however, seems to be less accurate than methods we
mentioned earlier, and it exhibits instability on higher spatial frequencies.

In the present paper we propose a novel non-iterative algorithm for the
numerical inversion of the spherical mean Radon transform from open surfaces
in 2-D. Our numerical experiments show that this method yields very accurate
reconstructions in the absence of noise, and is almost as stable as the
classical filtration/backprojection (FBP) algorithm (widely used for the
inversion of the regular Radon transform \cite{Kak,Natterer}). The present
algorithm requires pre-computation of certain functions that jointly serve as
a numerical filter; this needs to be done only once for a given configuration
of detectors. After these functions have been computed, each reconstruction
has numerical complexity similar to that of the FBP.

\section{Formulation of the problem}

In the thermoacoustic tomography acoustic detectors measure the pressure of
the outgoing wave radiating from the patient's body. This acoustic wave is
generated by the thermoacoustic expansion of the tissues initiated by a very
short electromagnetic pulse. The initial pressure $f(x)$ strongly depends on
the type of tissue, and is significantly higher for tumors, since they happen
to absorb much more electromagnetic energy than healthy tissue. Thus,
recovering $f(x)$ would yield important medical information about the location
and shape of tumors.

Let us denote by $g(z,t)$ the pressure registered at the moment $t$ by the
acoustic detector placed at the point $z$. Under certain assumptions (such as
a constant sound speed in the body, ideal infinitely small detector,
infinitely short pulse, and so on), one can recover from the measurements the
integral of $f(x)$ over a sphere of radius $r=ct$ centered at $z.$ A
dimensional analysis shows that in order to reconstruct a function of a 3-D
variable, the detectors should cover some measurement surface $\Sigma.$
Assuming for simplicity that the speed of sound equals unity, one arrives at
the following formulation of the inverse problem: reconstruct function $f(x)$
supported within some 3-D region of interest $\Omega$ from known values of its
integrals $g(z,r)$ over all spheres of radius $r$ with centers $z$ lying on
surface $S$ ($S$ and $\Omega$ are disjoint sets).

A 2-D version of this problem also arises in the thermoacoustic tomography.
Recently, it has been proposed to use integrating linear detectors instead of
point-like transducers \cite{Burghline,Paltaufline,Burghalg,Paltaufalg}. One
such detector consists of a long straight segment of optical fiber that serves
as a sensor of an optical interferometer; the measurements are proportional to
the integral the acoustic pressure over the length of the fiber. According to
the experimentalists, such detectors have much better sensitivity and spatial
resolution than the conventional transducer. It has been shown
\cite{Burghline} that the reconstruction problem in this case reduces to the
inversion of a set of the circular mean Radon transforms followed by the
inversion of a set of the regular 2-D Radon transforms (where the circular
mean Radon transform is a set of normalized integrals over circles with
centers lying on a curve) The same practical reasons as in 3-D case lead to
the requirement that the centers of the integration circles lie on an open
bounded curve satisfying (for a given ROI) the "visibility" condition.

In the present paper we study mostly the 2-D case; in order to simplify the
analysis we will concentrate on a truncated circular geometry as described
below. The 2-D region of interest $\Omega$ is a truncated disk $\Omega
(R,x_{\mathit{right}})=\{x=(x_{1},x_{2})\left\vert {}\right.  x_{1}^{2}%
+x_{2}^{2}<R^{2}\mathrm{\ and\ }x_{1}<x_{\mathit{right}}\}$ and the centers of
integration circles lie on the circular arch $\gamma(R_{\gamma}%
,z_{\mathit{right}})=\{z=(z_{1},z_{2})\left\vert {}\right.  z_{1}^{2}%
+z_{2}^{2}=R_{\gamma}^{2}\mathrm{\ and\ }z_{1}<z_{\mathit{right}}\},$ where
$R_{\gamma}>R,$ see Figure \ref{figgeom}. Our goal is to reconstruct a
$C_{0}^{1}$ function $f(x)$ supported in $\Omega$ from the known values of its
integrals $g(z,r)$ over circles $\mathbb{S}(r,z)$ of radii $r$ centered at
points $z\in\gamma$%
\[
g(z,r)=\int\limits_{\mathbb{S}(r,z)}f(x)dl(x)=r\int\limits_{\mathbb{S}^{1}%
}f(z+r\varpi)d\varpi.
\]
(The circular mean Radon transform is defined by the normalized integrals (or
means), i.e. by values of $g(z,r)/2\pi r.$ We, however, prefer to work with
the integrals $g(z,r)$ ). The invertibility of the circular mean Radon
transform is well known for such geometry \cite{AQ}. The stability of the
reconstruction problem is, again, determined by the "visibility" condition
(see \cite{LQ,XWAK1,XWAK2,Palambook,Q90}), which for this geometry is
equivalent to the inequality $z_{\mathit{right}}\geq x_{\mathit{right}}.$ Our
goal is to develop an efficient computational algorithm for the solution of
this problem.

\begin{figure}[t]
\begin{center}
\includegraphics[width=1.76in,height=2.2in]{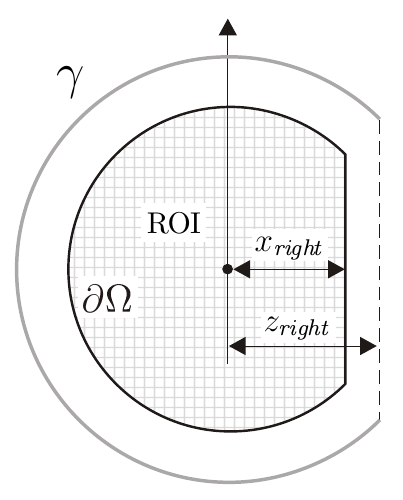}
\end{center}
\caption{Geometry}%
\end{figure}\label{figgeom}

\section{Outline of the method\label{method}}

The present algorithm is based on precomputing the approximations of plane
waves in $\Omega$ by the single layer potentials in the form $\int_{\gamma
}Z(\lambda|z-x|)\rho(z)dl(z),$ where $\rho(z)$ is the density of the potential
and $Z(t)\ $is either the Bessel function $J_{0}(t)$ or the Neumann
function\ $Y_{0}(t).$ Specifically, given a wavevector $\xi$ we find
numerically the densities $\rho_{\xi,J}(z)$ and $\rho_{\xi,Y}(z)$ of the
potentials
\begin{align}
W_{J}(x,\rho_{\xi,J}) &  =\int_{\gamma}J_{0}(\lambda|z-x|)\rho_{\xi
,J}(z)dl(z),\label{potj}\\
W_{Y}(x,\rho_{\xi,Y}) &  =\int_{\gamma}Y_{0}(\lambda|z-x|)\rho_{\xi
,Y}(z)dl(z),\label{poty}%
\end{align}
where $\lambda=|\xi|,$ such that
\begin{equation}
W_{J}(x,\rho_{\xi,J})+W_{Y}(x,\rho_{\xi,Y})\thickapprox\exp(-i\xi\cdot
x),\quad\forall x\in\Omega,\label{approx}%
\end{equation}
Obtaining such approximations is not trivial; we discuss this issue in more
detail in the next section. However, if the densities $\rho_{\xi,J}$ and
$\rho_{\xi,Y}$ have been found for all $\xi$ then function $f(x)$ can be
easily reconstructed. Indeed, let us introduce convolutions $G_{J}%
(\lambda,y),$ $G_{Y}(\lambda,z)$ as follows
\begin{align*}
G_{J}(\lambda,z) &  =\int_{\Omega}f(x)J_{0}(\lambda|z-x|)dx,\\
G_{Y}(\lambda,z) &  =\int_{\Omega}f(x)Y_{0}(\lambda|z-x|)dx.
\end{align*}
We notice that the boundary values of these functions for all $z\in\gamma$ can
be computed from projections $g(z,r)$:%
\begin{align}
G_{J}(\lambda,z) &  =\int_{R^{+}}g(z,r)J_{0}(\lambda r)dr,\label{gj}\\
G_{Y}(\lambda,z) &  =\int_{R^{+}}g(z,r)Y_{0}(\lambda r)dr.\label{gy}%
\end{align}
Consider now the Fourier transform $\hat{f}(\xi)$ of $f(x)$%
\[
\hat{f}(\xi)=\frac{1}{2\pi}\int_{\Omega}f(x)\exp(-i\xi\cdot x)dx.
\]
Using (\ref{approx}) we obtain
\begin{align}
\hat{f}(\xi) &  \thickapprox\frac{1}{2\pi}\int_{\Omega}f(x)\left[
W_{J}(x,\rho_{\xi,J})+W_{Y}(x,\rho_{\xi,Y})\right]  dx\nonumber\\
&  =\frac{1}{2\pi}\int_{\gamma}\left[  \int_{\Omega}f(x)J_{0}(\lambda
|z-x|)dx\right]  \rho_{\xi,J}(z)dl(z)\nonumber\\
&  +\frac{1}{2\pi}\int_{\gamma}\left[  \int_{\Omega}f(x)Y_{0}(\lambda
|z-x|)dx\right]  \rho_{\xi,Y}(z)dl(z)\nonumber\\
&  =\frac{1}{2\pi}\int_{\gamma}\left[  \rho_{\xi,J}(z)G_{J}(\lambda
,z)+\rho_{\xi,Y}(z)G_{Y}(\lambda,z)\right]  dl(z).\label{onefreq}%
\end{align}
Formulas (\ref{gj}) and (\ref{gy}) in combination with (\ref{onefreq}) allow
us to reconstruct (approximately) the Fourier transform $\hat{f}(\xi)$ of
$f(x).$ Now $f(x)$ can be recovered by inverting the 2D Fourier transform. In
order to minimize the operation count we will compute $G_{J}(\lambda,z)$ and
$G_{Y}(\lambda,z)$ (using (\ref{gj}) and (\ref{gy})) for a set of fixed values
of $\lambda=|\xi|;$ thus, the values of $\hat{f}(\xi)$ will be found on a
polar grid in $\xi.$ One can now use a high-order 2-D interpolation to obtain
values of the Fourier transform on the Cartesian grid and apply the inverse 2D
Fast Fourier Transform (FFT). Another approach is to utilize the famous
slice-projection theorem (see, for example \cite{Natterer}) and reconstruct
from the values of $\hat{f}(\xi)$ on the polar grid the regular Radon
projections of $f(x);$ the function then can be obtained by the application of
the FBP algorithm \cite{Kak,Natterer}. We have chosen the latter approach. The
whole algorithm can be briefly outlined as follows.

\begin{enumerate}
\item Choose a uniform polar grid in $\xi.$ For each value of $\xi
_{i,j}=\lambda_{i}(\cos\theta_{j},\sin\theta_{j})$ precompute $\rho_{\xi
_{i,j},J}$ and $\rho_{\xi_{i,j},Y}$ as described in the next section. This
step does not depend on values of $g(z,r)$ and needs to be performed only once
for each particular geometry.

\item Given values of $g(z,r),$ compute $G_{J}(\lambda_{i},z)$ and
$G_{Y}(\lambda_{i},z)$ (using (\ref{gj}) and (\ref{gy})) for each value of
$i.$ This can be done using the trapezoid rule in $r.$

\item Compute $\hat{f}(\xi_{i,j})$ according to (\ref{onefreq}), using the
same discretization of the integral as was used to obtain (\ref{potj}) and
(\ref{poty}).

\item Apply 1-D FFT to values of $\hat{f}(\xi_{i,j})$ for each fixed $j,$ to
reconstruct the classical Radon projections of $f(x)$ corresponding to the
angular parameter $\theta_{j}.$ (Alternatively, in the presence of strong
noise in the data one can introduce a lower-pass filter $\eta(\xi)$ and use
$\eta(|\xi_{i,j}|)\hat{f}(\xi_{i,j})$ instead of $\hat{f}(\xi_{i,j})$ on this step).

\item Use the well-know FBP algorithm \cite{Kak,Natterer} to reconstruct
$f(x)$ from the standard Radon projections.
\end{enumerate}

If we assume, for simplicity, that the size of the reconstruction grid is
$n\times n,$ and that the number of detectors and the number of integrals in
one projection are of the same order (say, $\mathcal{O}(n)),$ then steps 2, 3,
and 5 of the algorithm require $\mathcal{O}(n^{3})$ floating point operations;
step 5 is faster ($\mathcal{O}(n^{2}\log n)$ flops). Step 1, described in the
next section, is more expensive computationally. The precise operation count
depends on the algorithm used to compute $n$ Singular Value Decompositions
(SVD) for matrices of size $(n\times n),$ and can be $\mathcal{O}(n^{4})$
operations or higher. However, step 1 needs to be performed only once for each
particular geometry. The number of data generated on step 1 and stored on a
hard drive (values of the densities $\rho_{\xi_{i,j},J}$ and $\rho_{\xi
_{i,j},Y}$) is of order $\mathcal{O}(n^{3}).$ The rest of the algorithm,
including reading the precomputed densities, is completed in $\mathcal{O}%
(n^{3})$ operations (similarly to the FBP).

Obviously, the feasibility of this method hinges on our ability to find
approximations in the form (\ref{approx}). In the following section we discuss
the existence, accuracy and stability of such approximations.

\section{Approximations of plane waves by the single layer potentials}

\subsection{Full circle acquisition}

Before studying the more interesting case of an open acquisition curve
$\gamma$ let us consider a simpler case of the full-circle acquisition
geometry. Suppose the 2-D region of interest $\Omega$ is the disk of radius
$R$ centered at the origin and the detectors lie on the concentric circle
$\gamma$ of radius $R_{\gamma}>R$ (this formally corresponds to a particular
case when $x_{\mathit{right}}\geq R$ and $z_{\mathit{right}}\geq R_{\gamma}).$
In order to represent the plane wave $\exp(i\xi\cdot x)$ by the single layer
potentials supported on $\gamma$, we expand the wave in the Fourier series in
polar angle $\theta_{x}$ as follows:%
\begin{equation}
\exp(i\xi\cdot x)=\sum_{n=-\infty}^{+\infty}i^{n}J_{|n|}(\lambda
|x|)\exp(in\left[  \theta_{x}-\theta_{\xi}\right]  )) \label{Anger}%
\end{equation}
where $\xi=\lambda(\cos\theta_{\xi},\sin\theta_{\xi}),$ $x=|x|(\cos\theta
_{x},\sin\theta_{x})$ (this is a well known Jacobi-Anger
expansion\ \cite{Colton,Watson}). We also utilize the addition theorem
\cite{Colton,Watson} for the Hankel function $H_{0}^{(1)}(\cdot)=J_{0}%
(\cdot)+iY_{0}(\cdot)$:%
\begin{equation}
H_{0}^{(1)}(\lambda|x-z|)=\sum_{n=-\infty}^{+\infty}H_{|n|}^{(1)}%
(\lambda|z|)J_{|n|}(\lambda|x|)\exp(in\left[  \theta_{x}-\theta_{z}\right]  ),
\label{addition}%
\end{equation}
where $z=|z|(\cos\theta_{z},\sin\theta_{z}).$ By integrating equation
(\ref{addition}) with the density $\exp(in\theta_{z})$ supported on $\gamma$
one obtains
\[
\int_{\gamma}H_{0}^{(1)}(\lambda|z-x|)\exp(in\theta_{z})dl(z)=2\pi R_{\gamma
}H_{|n|}^{(1)}(\lambda R_{\gamma})J_{|n|}(\lambda|x|)\exp(in\theta_{x})
\]
or%
\[
\overline{H_{|n|}^{(1)}(\lambda R_{\gamma})}\int_{\gamma}H_{0}^{(1)}%
(\lambda|z-x|)\exp(in\theta_{z})dl(z)=2\pi R_{\gamma}\left\vert H_{|n|}%
^{(1)}(\lambda R_{\gamma})\right\vert ^{2}J_{|n|}(\lambda|x|)\exp(in\theta
_{x}).
\]
The latter expression can be used to express the cylindrical wave
$J_{|n|}(\lambda|x|)\exp(in\theta_{x})$ in the form%
\begin{align*}
J_{|n|}(\lambda|x|)\exp(in\theta_{x})  &  =\frac{J_{|n|}(\lambda R_{\gamma}%
)}{2\pi R_{\gamma}\left\vert H_{|n|}^{(1)}(\lambda R_{\gamma})\right\vert
^{2}}\int_{\gamma}J_{0}(\lambda|z-x|)\exp(in\theta_{z})dl(z)\\
&  +\frac{Y_{|n|}(\lambda R_{\gamma})}{2\pi R_{\gamma}\left\vert H_{|n|}%
^{(1)}(\lambda R_{\gamma})\right\vert ^{2}}\int_{\gamma}Y_{0}(\lambda
|z-x|)\exp(in\theta_{z})dl(z).
\end{align*}
Further, by substituting this formula into Jacobi-Anger expansion
(\ref{Anger}), one can represent $\exp(i\xi\cdot x)$ as the sum of the single
layer potentials%
\begin{align}
\exp(i\xi\cdot x)  &  =\frac{1}{2\pi R_{\gamma}}\int_{\gamma}\left[
\sum_{n=-\infty}^{+\infty}\frac{J_{|n|}(\lambda R_{\gamma})}{\left\vert
H_{|n|}^{(1)}(\lambda R_{\gamma})\right\vert ^{2}}i^{n}\exp(in\left[
\theta_{z}-\theta_{\xi}\right]  ))\right]  J_{0}(\lambda|z-x|)dl(z)\nonumber\\
&  +\frac{1}{2\pi R_{\gamma}}\int_{\gamma}\left[  \sum_{n=-\infty}^{+\infty
}\frac{Y_{|n|}(\lambda R_{\gamma})}{\left\vert H_{|n|}^{(1)}(\lambda
R_{\gamma})\right\vert ^{2}}i^{n}\exp(in\left[  \theta_{z}-\theta_{\xi
}\right]  ))\right]  Y_{0}(\lambda|z-x|)dl(z)\nonumber\\
&  =W_{J}(x,\rho_{\xi,J})+W_{Y}(x,\rho_{\xi,Y}) \label{exact}%
\end{align}
with the densities $\rho_{\xi,J}(z)$ and $\rho_{\xi,Y}(z)$ defined by the
formulas%
\begin{align*}
\rho_{\xi,J}(z)  &  =\frac{1}{2\pi R_{\gamma}}\sum_{n=-\infty}^{+\infty}%
\frac{J_{|n|}(\lambda R_{\gamma})}{\left\vert H_{|n|}^{(1)}(\lambda R_{\gamma
})\right\vert ^{2}}i^{n}\exp(in\left[  \theta_{z}-\theta_{\xi}\right]  )),\\
\rho_{\xi,Y}(z)  &  =\frac{1}{2\pi R_{\gamma}}\sum_{n=-\infty}^{+\infty}%
\frac{Y_{|n|}(\lambda R_{\gamma})}{\left\vert H_{|n|}^{(1)}(\lambda R_{\gamma
})\right\vert ^{2}}i^{n}\exp(in\left[  \theta_{z}-\theta_{\xi}\right]  )).
\end{align*}
The above series converge uniformly due to the fast growth of the Hankel
functions (see, for example \cite{Colton}) as the order $n$ goes to infinity.
Interestingly, in this simple case of the full-circle acquisition,
representation (\ref{exact}) is exact, and therefore the reconstruction
algorithm described in the previous section is also theoretically exact (in
this particular case).

For future use we compute the $L^{2}$ norm $N(\lambda)$ of the pair of the
densities $\left(  \rho_{\xi,J},\rho_{\xi,Y}\right)  $ defined as follows
\[
N^{2}(\lambda)\equiv||\left(  \rho_{\xi,J},\rho_{\xi,Y}\right)  ||_{2}%
^{2}=\int_{\gamma}\left[  \rho_{\xi,J}^{2}(z)+\rho_{\xi,Y}^{2}(z)\right]
dl(z).
\]
Due to the orthogonality of the complex exponents one obtains the following
simple formula%
\begin{equation}
N^{2}(\lambda)=\sum_{n=-\infty}^{+\infty}\frac{1}{\left\vert H_{|n|}%
^{(1)}(\lambda R_{\gamma})\right\vert ^{2}}. \label{treshold}%
\end{equation}
The values of $N(\lambda)$ can be easily computed numerically; it turns out
that for large values of $\lambda$ this function grows asymptotically linearly
and $N(\lambda)\approx\sqrt{2}\lambda.$

\subsection{Open curve $\gamma$: general considerations\label{general}}

The case when detectors are placed on an open curve $\gamma$ is more
complicated. In particular, in this case equation (\ref{approx}) cannot be
made into exact equality$.$ Indeed, the single layer potentials defined by
equations (\ref{potj}), (\ref{poty}) are solutions of the Helmholtz equations
in $\mathbb{R}^{2}\backslash\gamma;$ these functions decay at infinity as
$\mathcal{O}(|x|^{-1/2})$. The plane wave with $|\xi|=\lambda$ also solves the
same Helmholtz equation, but it does not decay. If these two solutions were
made to coincide within the open set $\Omega,$ they would also have the same
behavior at infinity, which is clearly impossible.

Not everything is lost, however. Let us consider the problem of approximating
$s(x)\in L^{2}(\partial\Omega)$ by a single layer potential $W_{H^{(1)}%
}(x,\rho)$ in the form
\begin{equation}
W_{H^{(1)}}(x,\rho)=\int_{\gamma}H_{0}^{(1)}(\lambda|z-x|)\rho(z)dl(z),
\label{hpotent}%
\end{equation}
where Hankel function $H_{0}^{(1)}(\lambda|z-x|)$ coincides (up to a constant
factor) with the free-space Green's function of the Helmholtz equation%
\[
\Delta u+\lambda^{2}u=0,
\]
subject to the radiation conditions at infinity. Assume additionally, that
$\lambda$ is not an eigenvalue of the Laplacian on $\Omega$ with zero boundary
conditions. If $s(x)$ cannot be approximated by potentials (\ref{hpotent}) in
the $L^{2}$ sense then there exists a non-zero function $t(x)\in
L^{2}(\partial\Omega)$ orthogonal to all such potentials:%
\[
\int_{\partial\Omega}t(x)\left[  \int_{\gamma}H_{0}^{(1)}(\lambda
|z-x|)\rho(z)dl(z)\right]  dl(x)=0,\qquad\forall\rho(z)\in L^{2}(\gamma).
\]
By interchanging the order of integration we obtain%
\[
\int_{\gamma}\rho(z)\left[  \int_{\partial\Omega}t(x)H_{0}^{(1)}%
(\lambda|z-x|)dl(x)\right]  dl(z)=0,\qquad\forall\rho(z)\in L^{2}(\gamma),
\]
which, in turn, implies%
\[
T(z)=\int_{\partial\Omega}t(x)H_{0}^{(1)}(\lambda|z-x|)dl(x)=0\qquad\forall
z\in\gamma.
\]
Function $T(z)$ defined by the above expression, is also a single layer
potential with density $t(x)$ supported on $\partial\Omega.$ This function is
real-analytic in $\mathbb{R}^{2}\backslash\bar{\Omega}.$ Since it vanishes on
$\gamma$ it must vanish on the whole circle $\mathbb{S}(r,0)=\left\{
x\left\vert {}\right.  x_{1}^{2}+x_{2}^{2}=R_{\gamma}^{2}\right\}  .$Thus,
$T(z)$ is the unique solution of the Dirichlet problem for the Helmholtz
equation \cite{Colton} in the exterior of the disk $\mathbb{B}(r,0)=\left\{
x\left\vert {}\right.  x_{1}^{2}+x_{2}^{2}<R_{\gamma}^{2}\right\}  $ circle
satisfying zero boundary conditions on $\mathbb{S}(r,0)$ and the radiation
condition at infinity. Therefore, $T(z)$ identically vanishes in
$\mathbb{R}^{2}\backslash\mathbb{B}(r,0)$. By the analyticity of the solutions
of the Helmholtz equation, $T(z)$ vanishes in $\mathbb{R}^{2}\backslash
\Omega.$ On the other hand, by the continuity of the single layer potentials,
$T(z)$ also approaches zero when $z$ and approaches $\partial\Omega$ from
inside of $\Omega$. Since $\lambda$ is not an eigenvalue of the Laplacian on
$\Omega,$ $T(z)$ vanishes in $\Omega.$ Now the well-known jump condition for
single layer potentials \cite{Colton} implies that $t(x)$ identically equals
zero on $\partial\Omega.$ This contradiction proves the following

\begin{theorem}
An arbitrary $s(x)\in L^{2}(\partial\Omega)$ can be approximated by potentials
in the form (\ref{hpotent}) in the $L^{2}$ sense. (In other words, single
layer potentials (\ref{hpotent}) are dense in $L^{2}(\partial\Omega).\}$
\end{theorem}

A very similar statement can be proven if one replaces Hankel's function
$H_{0}^{(1)}$ by $H_{0}^{(2)}$ in the above proof. Therefore, the sums of
potentials%
\[
W_{H^{(1)}}(x,\rho_{1})+W_{H^{(2)}}(x,\rho_{2})=\int_{\gamma}H_{0}%
^{(1)}(\lambda|z-x|)\rho_{1}(z)dl(z)+\int_{\gamma}H_{0}^{(2)}(\lambda
|z-x|)\rho_{2}(z)dl(z)
\]
with arbitrary densities $\rho_{1}(z)$ and $\rho_{2}(z)$ are also dense in
$L^{2}(\partial\Omega).$ Alternatively, due to the equations%
\[
H_{0}^{(1,2)}(t)=J_{0}(t)\pm iY_{0}(t),
\]
the sum $W_{J}(x,\rho_{J})+W_{Y}(x,\rho_{Y})$ can be used for approximation,
where
\begin{align*}
W_{J}(x,\rho_{J})  &  =\int_{\gamma}J_{0}(\lambda|z-x|)\rho_{J}(z)dl(z),\\
W_{Y}(x,\rho_{Y})  &  =\int_{\gamma}Y_{0}(\lambda|z-x|)\rho_{Y}(z)dl(z),
\end{align*}
again, under the assumption that $\lambda$ is not an eigenvalue of the
Dirichlet Laplacian.

Our goal is to approximate the plane waves in $\Omega$ by the single layer
potentials. Due to the uniqueness and stability of the Dirichlet problem for
the Helmholtz equation, it is enough to approximate the boundary values of
these functions on $\partial\Omega,$ which can be done as explained above ---
except for the case when $\lambda$ is the eigenvalue of the Dirichlet
Laplacian. In the latter case, one can approximate the normal derivative of
the target function by the normal derivative of the single layer potential. A
derivation similar to the one in the beginning of this section shows that the
normal derivatives of the layer potentials $\frac{\partial}{\partial
n(x)}\left(  W_{J}(x,\rho_{J})+W_{Y}(x,\rho_{Y})\right)  $ also form a dense
set in $L^{2}(\partial\Omega)$, if $\lambda$ is not an eigenvalue of the
Neumann Laplacian on $\Omega.$

Finally, since the eigenvalues of Laplacian on $\Omega$ may not be known in
advance, it is advantageous to use simultaneously both Neumann and Dirichlet
data. Namely, in order to approximate a solution $u(x)$ of the Helmholtz
equation on $\Omega$ by the single layer potentials, we form a vector function
$(u(x),\frac{\partial}{\partial n(x)}u(x))\rule[-4pt]{0.5pt}{15pt}%
_{\partial\Omega}$ and try to approximate it in the $L^{2}$ sense by the
functions in the form $\left(  W_{J}(x,\rho_{J})+W_{Y}(x,\rho_{Y}%
),\frac{\partial}{\partial n(x)}\left[  W_{J}(x,\rho_{J})+W_{Y}(x,\rho
_{Y})\right]  \right)  \rule[-7pt]{0.5pt}{18pt}_{\partial\Omega}$.

\subsection{Stability and regularization\label{sectstab}}

As established in the previous section, given a plane wave $\exp(-i\xi\cdot
x),$ one can can find a sequence of densities $\left(  \rho_{\xi,J}^{(k)}%
,\rho_{\xi,Y}^{(k)}\right)  $, $k=1,2,...$, such that the sum of potentials
$W_{J}(x,\rho_{\xi,J}^{(k)})+W_{Y}(x,\rho_{\xi,Y}^{(k)})$ converges to the
wave as described. However, in the case of open $\gamma$ the sequence of
densities themselves cannot have a limit in $L^{2}$ sense; if this limit
existed, the sum of potentials would exactly equal the plane wave, which
cannot happen, as discussed in the beginning of section \ref{general}.

Moreover, the sequence $\left(  \rho_{\xi,J}^{(k)},\rho_{\xi,Y}^{(k)}\right)
$ can start growing very fast in the $L^{2}$ norm defined by the formula%
\[
||\left(  \rho_{J},\rho_{Y}\right)  ||_{2}^{2}=\int_{\gamma}\left[  \left\vert
\rho_{J}(z)\right\vert ^{2}+\left\vert \rho_{Y}(z)\right\vert ^{2}\right]
dl(z)
\]
This, in particular, should occur in the situation when curve $\gamma$ does
not satisfy the visibility condition. Indeed, it is known that the
reconstruction problem is strongly unstable, if this condition is not
satisfied \cite{LQ,XWAK1,XWAK2,Palambook,Q90}. On the other hand, possible
instabilities in the present algorithm are associated with densities that are
strongly oscillating and large in the $L^{2}$ norm. Observe that in equation
(\ref{onefreq}) we compute the inner product of densities $\rho_{\xi,J},$
$\rho_{\xi,Y}$ with functions $G_{J}(\lambda,z)$ obtained from the measured
data. Any component of noise well correlated with $\rho_{\xi,J}$ and
$\rho_{\xi,Y}$ will be strongly amplified if these functions are large (in
$L^{2}$ norm). This amplification is the only source of instability in
reconstructing the value of $\hat{f}(\xi)$ at a particular point $\xi$ by our
method, and it has to occur to render the reconstruction unstable --- as it
should be in accordance with the theory.

Our goal, of course, is to use the present technique with the data acquisition
configurations that satisfy the visibility conditions. It is difficult to find
theoretically a sharp estimate on the behavior of some "optimal" sequence of
densities that would combine convergence of potentials with the slow or no
growth of the densities. Instead, we propose a regularized algorithm for
computation of "good" approximations.

Let us introduce the families $\Gamma$ and $\Upsilon$ of pairs of $L^{2}$
functions defined on $\gamma$ and $\partial\Omega$:%
\begin{align*}
\Gamma &  =\{(q_{1}(z),q_{2}(z))|q_{1},q_{2}\in L^{2}(\gamma)\},\\
\Upsilon &  =\{(p_{1}(z),p_{2}(z))|p_{1},p_{2}\in L^{2}(\partial\Omega)\}.
\end{align*}
Define the inner products for functions from $\Gamma$ and $\Upsilon$ as
follows%
\begin{align*}
\forall\mathbf{q,s}  &  \mathbf{\in}\Gamma,\qquad\left\langle \mathbf{q,s}%
\right\rangle _{\Gamma}=\int_{\gamma}\left[  q_{1}(z)\overline{s_{1}(z)}%
+q_{2}(z)\overline{s_{2}(z)}\right]  dl(z),\\
\forall\mathbf{p,r}  &  \mathbf{\in}\Upsilon,\qquad\left\langle \mathbf{p,r}%
\right\rangle _{\Upsilon}=\int_{_{\partial\Omega}.}\left[  p_{1}%
(x)\overline{r_{1}(x)}+p_{2}(x)\overline{r_{2}(x)}\right]  dl(x).
\end{align*}
Further, let us introduce the operator $\mathbf{A}$ that maps a pair
$\mathbf{q\in}\Gamma$ into a pair $\mathbf{p\in}\Upsilon$ according to the
following formula%
\[
\mathbf{p=A}(\mathbf{q})\equiv\left(  W_{J}(x,q_{1})+W_{Y}(x,q_{2}),\frac
{1}{\lambda}\frac{\partial}{\partial n(x)}\left[  W_{J}(x,q_{1})+W_{Y}%
(x,q_{2})\right]  \right)  \rule[-10pt]{0.5pt}{30pt}_{\partial\Omega},
\]
where single layer potentials $W_{J}(x,q_{1}),$ $W_{Y}(x,q_{2})$ are defined,
as before, by equations (\ref{potj}), (\ref{poty}).

Given boundary values of the plane wave $\mathbf{u}_{\xi}(x)=(\exp(i\xi\cdot
x),\frac{1}{\lambda}\frac{\partial}{\partial n}\exp(i\xi\cdot x))_{\partial
\Omega}$ we would like to find a pair of functions (densities) $\left(
\rho_{\xi,J}(z),\rho_{\xi,Y}(z)\right)  $ such that
\begin{equation}
\mathbf{A}(\left(  \rho_{\xi,J},\rho_{\xi,Y}\right)  )\approx\mathbf{u}_{\xi},
\label{svdeq}%
\end{equation}
subject to the requirement that the densities are "not too large". There is
more than one way to find such regularized solutions; we utilize the SVD\ of
the operator $\mathbf{A.}$ Namely, we find (numerically) the sets of pairs
(left and right singular vectors) $\mathbf{q}^{(j)}\in\Gamma,$ $\mathbf{p}%
^{(j)}\in\Upsilon,$ $j=1,2,..,$ and singular values $\sigma_{j},$ $\sigma
_{1}\geq\sigma_{2}\geq...\geq\sigma_{n}\geq...,$ such that%
\begin{align*}
\left\langle \mathbf{q}^{(i)},\mathbf{q}^{(j)}\right\rangle _{\Gamma}  &
=\delta_{i,j},\\
\left\langle \mathbf{p}^{(i)},\mathbf{p}^{(j)}\right\rangle _{\Upsilon}  &
=\delta_{i,j},\\
\mathbf{p}^{(j)}  &  \equiv\sigma_{j}\mathbf{A}(\mathbf{q}_{j}),
\end{align*}
where $\delta_{i,j}$ is the Kronecker symbol ($\delta_{i,j}=0$ if $i\neq j,$
and $\delta_{i,j}=1$ otherwise). Now the desired approximation is given by the sum%

\begin{equation}
\left(  \rho_{\xi,J}(z),\rho_{\xi,Y}(z)\right)  \approx\sum_{j=1}^{j_{\max}%
}\mathbf{q}^{(j)}(x)\frac{1}{\sigma_{j}}\left\langle \mathbf{u}_{\xi
},\mathbf{p}^{(j)}\right\rangle _{\Upsilon}, \label{inversesvd}%
\end{equation}
where $j_{\max}$ serves as a regularization parameter. The increase in
$j_{\max}$ yields a closer $L^{2}$ fit of the single layer potential to the
plane wave, but it also may (and in certain cases will) lead to the unbounded
growth of the densities.

A frequently used regularization technique for solving ill-posed problems
using the SVD decomposition is to drop from (\ref{inversesvd}) all the terms
with $\sigma_{j}$ smaller than a certain threshold $\sigma_{\min}$ (and thus
define the $j_{\max})$. However, it is not clear how to choose such
$\sigma_{\min}$ optimally. Moreover, in general $\sigma_{\min}$ should depend
on the frequency $\lambda,$ since the norm of the densities may grow with
$\lambda$ as we saw in the example of circular acquisition geometry.

Instead of chosing $\sigma_{\min}$, we propose to use the circular case as a
benchmark to determine the regularization parameter $j_{\max}$. We found in
the numerical experiments that very good results are obtained when $j_{\max}$
is chosen to be the largest number such that $||\left(  \rho_{\xi,J},\rho
_{\xi,Y}\right)  ||_{2}<KN(\lambda)$, where $N(\lambda)$ is computed using
equation (\ref{treshold}) for the circular case with the same values of $R$
and $R_{\gamma},$ and $K>1$ is a constant; in all numerical experiments
presented below $K$ was equal to $1.5$.

The most computationally expensive step in finding the densities is the SVD
decomposition; however, the same SVD is used for all plane waves with the same
frequency $\lambda=|\xi|.$ Algorithmically, we computed the SVD by first
discretizing the integrals that define operator $\mathbf{A,}$ and by enforcing
equation (\ref{svdeq}) at some set of collocation points on $\partial\Omega.$
This results in a matrix that represents a discretized version of $\mathbf{A}%
$; the SVD is then computed using subroutine DGESVD from LAPACK. In practice,
the discretization of the integrals is dictated by the number and location of
the detectors; we assumed that they were distributed uniformly over $\gamma.$
We also utilized equispaced collocation points on $\partial\Omega;$ the number
of points was chosen to be twice the number of the detectors.

One can notice that although this approach guarantees bounded solutions
$\left(  \rho_{\xi,J},\rho_{\xi,Y}\right)  ,$ there is no theoretical estimate
on the accuracy of the approximation of the plane waves by the corresponding
potentials. However, when the densities have been found, it is very easy to
compute numerically the approximation error, and, if necessary, to re-run the
computation with different values of parameters.

The computation of the densities $\left(  \rho_{\xi,J}(z),\rho_{\xi
,Y}(z)\right)  $ constitutes the first step of the reconstruction algorithm
(see section \ref{method}). It is rather time consuming; however, for a
particular acquisition geometry it needs to be done only once.

\section{Numerical examples}

\begin{figure}[t]
\begin{center}%
\begin{tabular}
[c]{cc}%
\includegraphics[width=2in,height=2in]{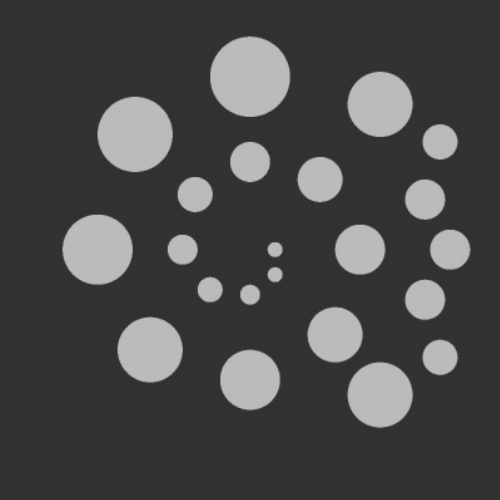} &
\includegraphics[width=2in,height=2in]{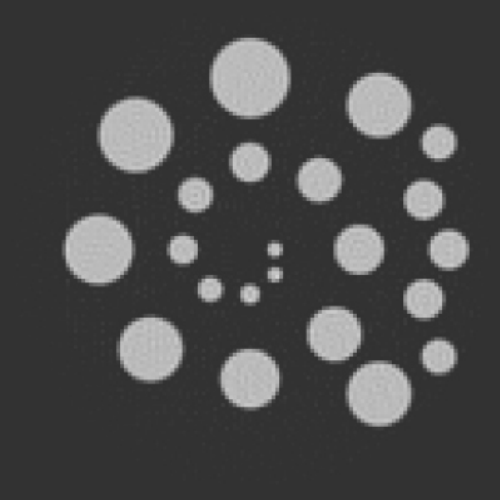}\\
(a) & (b)\\
& \\
\includegraphics[width=2in,height=2in]{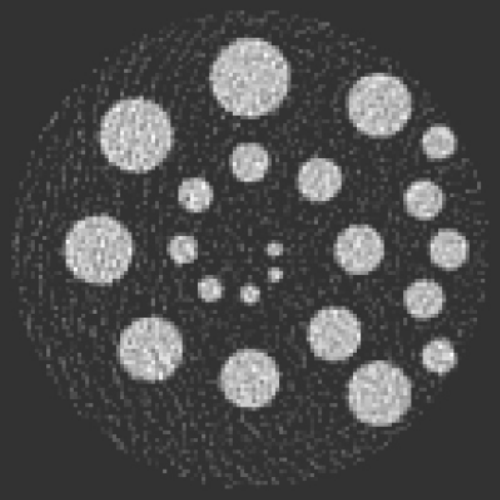} &
\includegraphics[width=2in,height=2in]{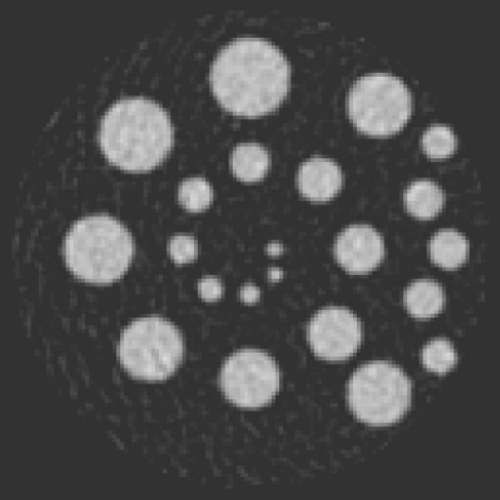}\\
(c) & (d)
\end{tabular}
\end{center}
\caption{Geometry \#1: (a) phantom, (b) reconstruction from the accurate
projections (c) effect of 15\% noise in projections (d) reconstruction from
the noisy projections with the additional filter}%
\label{figfull}%
\end{figure}

In order to verify the accuracy and stability of the present algorithm we
conducted a series of numerical experiments . In particular, we studied two
acquisition geometries (which we will call \#1 and \#2) corresponding to
different values of parameters $x_{\mathit{right}}$ and $z_{\mathit{right}},$
(see Figure \ref{figgeom}). In both geometries $R$ and $R_{\gamma}$ are equal
to $1$ and $1.3$ respectively. Geometry \#1 has $x_{\mathit{right}}%
$=$z_{\mathit{right}}=1;$ in other words, the ROI is a unit circle. Geometry
\#2 is defined by values $x_{\mathit{right}}$=$z_{\mathit{right}}=0,$ which
corresponds to a half-circle acquisition curve and a half-disk ROI. Both
geometries satisfy the visibility condition. In all experiments reconstruction
was performed on an equispaced Cartesian $129\times129$ grid, from simulated
projections corresponding to $500$ equispaced detectors, each measuring 129
circular integrals with equispaced radii. (The radial step in the projections
coincides with the step of the reconstruction grid under such discretization).

The goal of our first test was to verify the feasibility of the approximation
of plane waves by the single layer potentials. We found that the most
difficult plane waves to approximate were the ones propagating in the vertical
direction. In geometry \#1 the method described in Section \ref{sectstab} with
$K=1.5$ gave approximations accurate up to 4 decimal places or better. In
particular, for the plane wave propagating in the vertical direction with the
wavelength corresponding to the Nyquist frequency of our $129\times129$ grid,
the maximum pointwise error did not exceed $8\cdot10^{-6}.$ Similar (although
slightly less accurate) results were obtained in the geometry \#2.

These accurate approximations of the plane waves lead to very accurate
reconstructions of smooth images. For geometry \#1 we defined a smooth phantom
as follows. First, we introduced a $C_{0}^{8}(\mathbb{R})$ function $h(t)$:%
\[
h(t)=\left\{
\begin{array}
[c]{lll}%
c_{0}\int_{0}^{1-|t|}\sin^{8}(\pi s)ds & , & 0\leq|t|\leq1\\
0 & , & |t|>1
\end{array}
\right.
\]
where constant $c_{0}$ was chosen so that $h(0)=1.$ Function $h(t)$ is even,
compactly supported on $[-1,1],$ 8 times continuously differentiable on
$\mathbb{R}$ function, with $h(1/2)=1/2.$ We then used this function to
construct the 8 times continuously differentiable phantom $f(x)=h\left(
\frac{|x-x_{1}|}{r_{1}}\right)  +h\left(  \frac{|x-x_{2}|}{r_{2}}\right)  ,$
as a sum of two bell-shaped rotationally invariant functions. For geometry \#1
we used parameters $x_{1}=(0.3,0.3),x_{2}=(-0.4,0.2),$ $r_{1}=0.55,$
$r_{2}=0.5.$ Reconstruction using the present method (without additional
filtration) resulted in maximal point-wise error of $7.3\cdot10^{-5}.$
Similarly small reconstruction errors were obtained for smooth phantoms in
geometry \#2. We do not present here the gray scale pictures of these phantoms
and the corresponding reconstructed images since they look identically. The
experiments with smooth phantoms show that our method indeed reconstructs
accurately all the spatial frequencies of the image (including the lower ones)
--- as opposed to parametrix approximations.

Our next experiment was with discontinuous phantoms. Point-wise accurate
reconstructions are not possible with such functions due to the aliasing
errors and the Gibbs phenomenon; the goal is to obtain images that appear
qualitatively correct, with low noise amplification. As a phantom in geometry
\#1, we used the sum of the characteristic functions of circles as shown in
Figure \ref{figfull}(a). Part (b) of the latter figure demonstrates the image
reconstructed from accurate projections. In order to analyze the sensitivity
of the method to non-exact measurements we added to the projections white
noise with intensity of $15\%$ of the signal (in $L^{2}$ norm); the
reconstruction is shown in Figure \ref{figfull}(c). As it is frequently done
in tomography, in order to reduce the effects of noise one can apply a
low-pass filter on step 4 of the algorithm. Part (d) demonstrates the effect
of such additional filtration, with filter $\eta(\xi)=\cos(\frac{\pi}{2}%
|\xi|/\lambda_{\mathit{Nyquist}})$, where $\lambda_{\mathit{Nyquist}}$ is the
Nyquist frequency of the reconstruction grid.\begin{figure}[t]
\begin{center}%
\begin{tabular}
[c]{cc}%
\includegraphics[width=2in,height=2in]{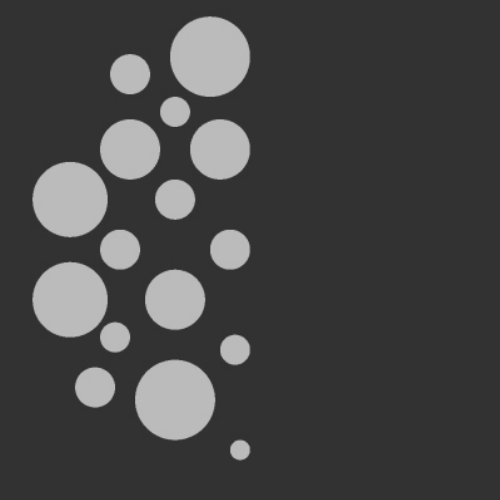} &
\includegraphics[width=2in,height=2in]{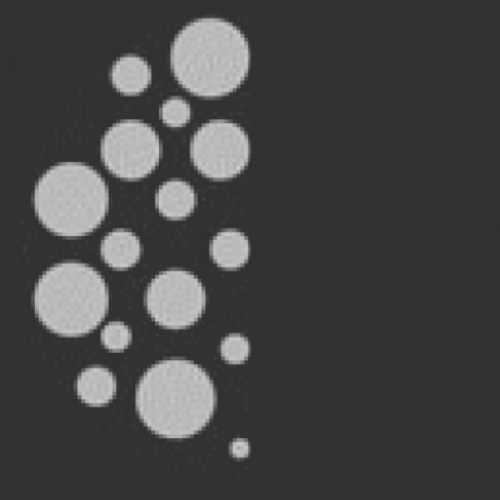}\\
(a) & (b)\\
& \\
\includegraphics[width=2in,height=2in]{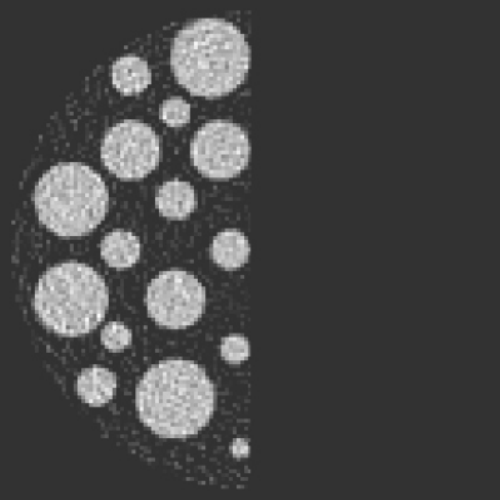} &
\includegraphics[width=2in,height=2in]{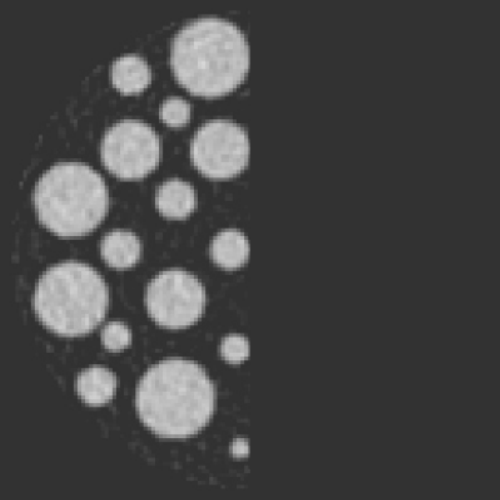}\\
(c) & (d)
\end{tabular}
\end{center}
\caption{Geometry \#2: (a) phantom, (b) reconstruction from the accurate
projections (c) effect of 15\% noise in projections (d) reconstruction from
the noisy projections with the additional filter}%
\label{fighalfb}%
\end{figure}


Quite similar results were obtained in geometry \#2. As a phantom we used
again a sum of the characteristic functions of circles supported within the
ROI (the left half of the unit disk)), as shown in Figure \ref{fighalfb}(a).
Parts (b), (c), and (d) of that figure demonstrate, respectively,
reconstruction results from the accurate projections, noisy projections, and
reconstruction from noisy projections with additional filtration (using the
same filter as in the previous example).

In order to compare stability of our algorithm to that of the classical FBP,
we have reconstructed the two phantoms (shown in Figure \ref{figfull}(a) and
Figure \ref{fighalfb}(a)) using the latter method, from a similar number of
the standard Radon projections with the same level of noise (not shown here).
The general quality and the level of noise in the reconstructions was quite
close to those of the images obtained by the present method.

\begin{figure}[t]
\begin{center}%
\begin{tabular}
[c]{cc}%
\includegraphics[width=2in,height=2in]{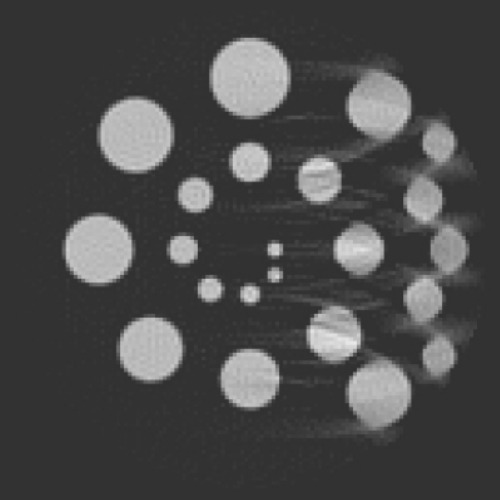} &
\includegraphics[width=2in,height=2in]{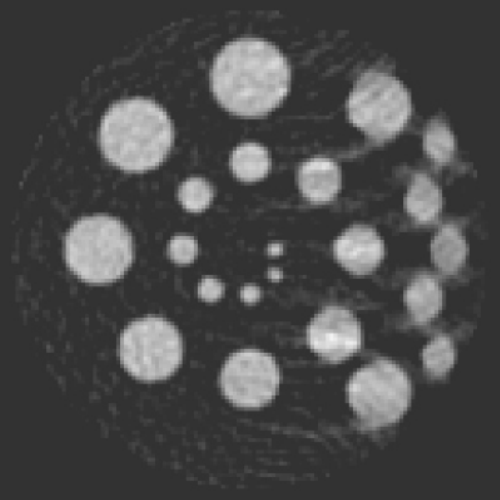}\\
(a) & (b)\\
&
\end{tabular}
\end{center}
\caption{Geometry \#2: (a) reconstruction with some sources outside of the ROI
(compare to Figure \ref{figfull}(a)) (b) the same with 15\% noise in the
projections and with the additional filter}%
\label{figbad}%
\end{figure}

From the practical point of view it is interesting to know how the
reconstruction is affected by acoustic sources located outside of the ROI. It
is known that, with the exception of the methods based on expansion in the
eigenvalues of the Dirichlet Laplacian on a closed domain
\cite{AKU,Kunyansky1}, all the other exact reconstruction techniques will
produce incorrect results in the presence of such sources (see
\cite{AKK,KuchKu,HKN} for further discussion of this phenomenon). In our case,
such sources will be present if the electromagnetic wave impinges on the parts
of the patient's body located outside the ROI. In other to model this
situation in geometry \#2 we used the previously used phantom supported within
the unit disk, as shown in Figure \ref{figfull}(a). The result of the
reconstruction in the absence of noise is shown in Figure \ref{figbad}(a);
part (b) demonstrates the effect of 15\% noise and the additional filtration.
One can notice in this figure that the parts of the source located
outside ROI are not reconstructed correctly (not unexpectedly). Nevertheless,
by covering the right half of the image, it is easy to see that the
reconstruction within the ROI is very little (if at all) affected by the
presence of additional sources outside the region. Interestingly, in the right
half of the image (outside ROI), the "visible" material interfaces are
reconstructed quite well, while the "invisible" ones (such that the normal
does not intersect the measuring surface) are noticeably smeared.

In accordance with the theoretical operation count, step \#1 turned out to be
the most expensive part of the algorithm; in the experiments described in this
section the computation of the densities took several hours. However, this
step needs to be performed only once for a given geometry. On consecutive runs
of the reconstruction program, once the pre-computed densities had been read
from the hard drive, the algorithm took a fraction of a second to complete.
The time required to read the densities was about 4 seconds. If our method
were to be used to process the data obtained by an acquisition scheme with
linear detectors (which requires the inversion of a set of 2-D problems), the
densities would need to be read only once, and each of the 2-D problems would
be inverted in a fraction of second.

\section{Concluding remarks}

We have presented an efficient reconstruction algorithm applicable to problems
of thermoacoustic tomography with the detectors lying on an open curve. The
method is based on $L^{2}$ approximations of plane waves by certain single
layer potentials; we have proven that such approximations are possible if the
measurement curve is an open circular arch. We have also verified numerically
that, for the truncated circular geometry satisfying the "visibility"
condition one can obtain very accurate approximations with bounded densities,
which, in turn, leads to stable image reconstructions.

In conclusion, we would like to add several remarks:

\begin{itemize}
\item The theorem we have presented does not restrict the shape of the ROI; it
can be arbitrary. Of course, for stable reconstruction the combination of the
ROI and the acquisition surface should satisfy the "visibility" condition.

\item In the proof of the theorem we used the fact that if the solution of the
Helmholtz equation vanishes on a circular arch, it must vanish on the whole
circle, and, therefore, in the exterior of the disk. A similar statement is
also true if the acquisition curve is a part of any closed analytic curve, and
the theorem extends to such configurations. The theorem also remains valid if
the measurement curve is a continuous non-analytic curve.

\item In order to use this approach with non-circular acquisition curves, the
regularization technique we presented may require some modifications, since it
is based on the comparison with the full-circle case.

\item The technique we propose can also be used in 3-D, for image
reconstruction from the detectors lying on an open surface. Compliance with
the "visibility" condition is, again, necessary for a stable reconstruction in
this case. We anticipate that the pre-computation of the densities may become
quite time-consuming in 3-D, due to the increased dimensionality of the
problem. The development of a fast algorithm for this part of our method is
required to make this approach practical (in 3-D). This will be the object of
our future research.
\end{itemize}

\section{Acknowledgements}

The author would like to thank Professor P. Kuchment for fruitful discussions
and numerous helpful suggestions.

This work was partially supported in part by the NSF/DMS grant NSF-0312292 and
by the DOE grant DE-FG02-03ER25577.

\end{document}